\documentclass[leqno]{article}
\usepackage{amsmath,amssymb}
\usepackage{latexsym}
         
    \def \Real{\mbox{\sl I\kern-.166em R}}            % real numbers
    \def \Hyp{\mbox{\sl I\kern-.166em H}}           
    \def \Nat{\mbox{\sl I\kern-.166em N}}             % natural numbers
    \newtheorem{theorem}{Theorem}[section]
    \newtheorem{proposition}[theorem]{Proposition}

    \newtheorem{definition}[theorem]{Definition}

    \newcommand{\red}{\protect\large\bf}

    \def\refprep#1#2#3{#1, #2, preprint, #3.}
    
    \def\refart#1#2#3#4{#1, #2, {\em #3\/}\ #4.}
    \def\book#1#2#3{#1, {\em #2,} #3.}

\title{\bf Homogeneous geodesics of non-unimodular Lorentzian Lie groups and naturally reductive Lorentzian spaces in dimension three}
\author{G. Calvaruso and  R.A. Marinosci%
\thanks{Authors supported by funds of MURST, GNSAGA and the University of Lecce.
\newline 2000 {\em Mathematics Subject Classification:} 53C20,53C22, 53C30.
\newline
{\em Keywords and phrases:} Lorentzian homogeneous spaces, homogeneous
geodesics, naturally reductive spaces, g.o. spaces. }
}
   \date{}

\begin{document}

    \maketitle

\begin{abstract}
We determine, for all three-dimensional non-unimodular Lie groups equipped with a Lorentzian metric, the set of homogeneous geodesics through a point. Together with the results of [C] and [CM2], this leads to the full classification of three-dimensional Lorentzian g.o. spaces and naturally reductive spaces. 
\end{abstract}

    \bigskip\noindent
    \section{\red Introduction}

A (connected) pseudo-Riemannian manifold $(M,g)$ is {\em homogeneous} provided that there exists a group $K$ of isometries acting transitively on it [O'N], that is, for any points $p,q \in M$, there is an isometry $\phi \in K$ such that $\phi (p)=q$. Such $(M,g)$ can be then identified with $(K/H,g)$, where $H$ is the isotropy group at a fixed point $o$ of $M$. 
We recall here a few results concerning homogeneous manifolds, in the Riemannian and pseudo-Riemannian case (in particular, in  Lorentzian geometry). 

Gadea and Oubi\~na [GO] introduced {\em homogeneous pseudo-Riemannian structures} in order to characterize reductive homogeneous pseudo-Riemannian manifolds. Earlier, a corresponding result had been proved for all homogeneous Riemannian manifolds by Ambrose and Singer [AS] (see also [TV]). 

Sekigawa [Se] proved that a three-dimensional connected, simply connected and complete homogeneous Riemannian manifold is either symmetric or it is a Lie group endowed of a left-invariant Riemannian metric. Taking into account the classification of three-dimensional Riemannian Lie groups given by Milnor [Mi], this result permits to determine all three-dimensional homogeneous Riemannian manifolds.      

Recently, the first author obtained the following similar result in the Lorentzian framework:

\begin{theorem}{\bf [C]}
Let $(M,g)$ be a three-dimensional connected, simply connected, complete homogeneous Lorentzian manifold. Then, either $(M,g)$ is symmetric, or it is isometric to a three-dimensional Lie group equipped with a left-invariant Lorentzian metric.
\end{theorem}

Theorem 1.1, together with the results on three-dimensional Lorentzian Lie groups obtained by Cordero and Parker [C0Pa2] and Rahmani [R], leads to the classification of three-dimensional homogeneous Lorentzian manifolds. 

It is then natural to try to characterize and classify some special classes of three-dimensional homogeneous Lorentzian manifolds having a special geometric meaning, also in order to compare such results in the Lorentzian case with their Riemannian analogues. In [C], the first author classified three-dimensional Lorentzian symmetric spaces.  Broader interesting classes of homogeneous Lorentzian spaces are given by naturally reductive and g.o. spaces, both related to the notion of homogeneous geodesic. 

Let $(M=K/H,g)$ be a pseudo-Riemannian homogeneous manifold. A geodesic $\Gamma$ through the origin $o \in M=K/H$ is called {\em homogeneous} if it is the orbit of a $1$-parameter subgroup. In general, the group $K$ is not unique. If $\Gamma$ is homogeneous with respect to some isometry group $K'$, then it is also homogeneous with respect to the maximal connected group of isometries $K$, but the converse does not hold [KNVl]. 

Homogeneous geodesics of homogeneous Riemannian manifolds have been investigated by many authors. We can refer to [CM1],[CKM],[K2],[KNVl],[KS],[KV],[KVl],[M], for some examples and further references. In [KS], O. Kowalski and J. Szenthe proved the fundamental result that  {\em any homogeneous Riemannian manifold admits at least one homogeneous
geodesic}. 

A corresponding existence result holds in the Lorentzian case, provided that the space is reductive homogeneous [P]. In the framework of Lorentzian geometry, homogeneous geodesics also aquire a new interest, since homogeneous Lorentzian spaces for which all null geodesics are homogeneous, are candidates for constructing solutions to the $11$-dimensional supergravity, which preserve more than $24$ of the available $32$ supersymmetries. In fact, all Penrose limits, preserving the amount of supersymmetry, of such a solution, must preserve homogeneity, which is the case for the Penrose limit of a reductive homogeneous spacetime along a null homogeneous geodesic ([Me],[P],[FMeP]).

A pseudo-Riemannian reductive homogeneous space $(M=K/H,g)$ is called a {\em g.o. space} if all its geodesics are homogeneous, it is  {\em naturally reductive} if there exists at least one reductive  split $\mathfrak{k}$ = $\mathfrak{m} \oplus \mathfrak{h}$ such that   
\begin {equation}
<[X,Y]_{\mathfrak{m}},  Z> + <[X,Z]_{\mathfrak{m}}, Y>=0,
\end{equation}
for all $X,Y,Z \in \mathfrak{m}$. It is not always easy to decide whether a homogeneous (reductive) pseudo-Riemannian manifold is or is not naturally reductive, because condition (1.1) must be checked for all groups of isometries acting transitively on $M$ [TV]. It is also well-known that (1.1) holds if and only if the Levi-Civita connection of $(M,g)$ and the canonical connection (of the reductive split $\mathfrak{k}$ = $\mathfrak{m} \oplus \mathfrak{h}$) have exactly the same geodesics [TV].  Clearly, a naturally reductive space is g.o., but in dimension greater than $5$ there exist Riemannian g.o. spaces which are in no way naturally reductive [KV], while in dimension three the two classes of Riemannian homogeneous spaces coincide. It is well-known that symmetric spaces are special examples of naturally reductive spaces.  We recall that three-dimensional naturally reductive Lorentzian spaces have been already investigated by Cordero and Parker in [CoPa1], in order to determine the possible forms and the symmetry groups of their curvature tensor.
 
In [CM2], the authors determined homogeneous geodesics of all three-dimensional unimodular Lie groups admitting a left-invariant Lorentzian metric. In this paper, we investigate  homogeneous geodesics of three-dimensional non-unimodular Lorentzian Lie groups. Taking into account Theorem 1.1, this permits to determine all three-dimensional g.o. and naturally reductive Lorentzian spaces.

The paper is organized in the following way. In Section 2, we shall recall the basic definitions and properties of homogeneous geodesics in a homogeneous pseudo-Riemannian manifold. In Section 3, we shall report the classification of three-dimensional homogeneous Lorentzian manifolds, and we shall describe the set of geodesic vectors for all three-dimensional non-unimodular Lorentzian Lie groups. In Section 4 we will then give the classification of three-dimensional naturally reductive and g.o. Lorentzian manifolds.

\bigskip \noindent
 \section{\red Preliminaries on homogeneous geodesics}
 \setcounter{equation}{0}

Let $(M,g)$ be a (connected) homogeneous pseudo-Riemannian manifold. Then, its full isometry group $I(M)$ acts transitively on it and $M$ can be identified with $(K/H,g)$, where $K \subset I(M)$ is a connected subgroup of $I(M)$ acting transitively on $M$ and $H$ is the isotropy group at a fixed point $o\in M$. In general, we can have more than one choice for $K$. 
 
In contrast to the Riemannian case, the Lie algebra $\mathfrak{k}$ of $K$ does not need to admit a reductive decomposition. Denote by $\mathfrak{k}$ and $\mathfrak{h}$ the Lie algebras of $K$ and $H$ respectively, and let $\mathfrak{m}$ be a complement of $\mathfrak{h}$ in $\mathfrak{k}$. If $\mathfrak{m}$ is stable under the action of $\mathfrak{h}$, then $\mathfrak{k}$ = $\mathfrak{m} \oplus \mathfrak{h}$ is called a {\em reductive split}, and $(\mathfrak{k}, \mathfrak{h})$ a {\em reductive pair}. It is important to stress that reductivity is not an intrinsic property of $(M,g)$, but of the description of $M$ as coset space $K/H$. In fact, the socalled {\em Kaigorodov space} is an example of a homogeneous Lorentzian manifold which has two different coset descriptions, but only one of them is reductive [FMeP]. Nevertheless, a homogeneous pseudo-Riemannian manifold $(M,g)$ is called reductive if there exists a Lie group $K$ acting transitively on $M$ via isometries, with corresponding isotropy group $H$, such that $(\mathfrak{k}, \mathfrak{h})$ is reductive. 

Reductive homogeneous pseudo-Riemannian manifolds are characterized by the existence of a pseudo-Riemannian homogeneous structure. 
Let $M$ be a connected manifold and $g$ a pseudo-Riemannian metric on $M$.
We denote by $\nabla$ the Levi-Civita connection of $(M,g)$ and by $R$ its curvature tensor. 

 \begin{definition}
{\em A} homogeneous pseudo-Riemannian structure {\em on $(M,g)$ is a tensor field $T$ of type $(1,2)$ on $M$,
such that the connection $\tilde{\nabla} = \nabla -T$ satisfies}
 \begin{equation}
\tilde{\nabla} g=0, \qquad \tilde{\nabla} R=0, \qquad \tilde{\nabla} T =0.
 \end{equation}
 
 \end{definition}
 
 In the Riemannian case, homogeneous structures were first introduced by Ambrose and Singer [AS], and further investigated by Tricerri and Vanhecke [TV]. Gadea and Oubi\~{n}a proved the following
 
 \begin{theorem}{\bf [GO]}
Let $(M,g)$ be a connected, simply connected and complete pseudo-Rieman- nian manifold. Then, $(M,g)$ admits a homogeneous pseudo-Riemannian structure if and only if it is a reductive homogeneous pseudo-Riemannian manifold.
 \end{theorem}
 
The connection $\tilde{\nabla}$, satisfying conditions (2.1), is called the {\em canonical connection} associated to the homogeneous pseudo-Riemannian structure $T$. In the special case of a symmetric space, the torsion tensor $\tilde{T}$ of $\tilde{\nabla}$ satisfies  the condition $\tilde{T}=0$ and hence, $\tilde{\nabla}$ coincides with the Levi-Civita connection $\nabla$ of $(M,g)$. In [C], in order to prove Theorem 1.1, the first author showed that any non-symmetric three-dimensional homogeneous Lorentzian manifold $(M,g)$ admits a homogeneous Lorentzian structure $T$ such that $T_X Y =\nabla _X Y$ for all $X,Y$ vector fields tangent to $M$. In particular, by Theorem 2.2, all three-dimensional homogeneous Lorentzian spaces are then reductive.

Let $\mathfrak{k}$ = $\mathfrak{m} \oplus \mathfrak{h}$ a reductive split coming from a homogeneous pseudo-Riemannian structure $T$. The geodesics of $\tilde{\nabla} = \nabla -T$ are curves of the form
\begin{equation}
\Gamma(t) = exp(tX)(o),\hspace{5pt}       t \in \Real ,
\end{equation}
with $X \in \mathfrak{m}$. They are called {\em canonical geodesics} of $(M,g)$ [KoNo,Chapter 10, Cor.2.5].

Consider now a reductive homogeneous pseudo-Riemannian manifold $(M=K/H,g)$, where 
$\mathfrak{k}$ = $\mathfrak{m} \oplus \mathfrak{h}$ is a reductive split. As already mentioned in the Introduction, a geodesic $\Gamma$ through the origin $o \in K/H$ is homogeneous if it is of the form $\Gamma(t) = exp(tZ)(o)$, for some $Z \in \mathfrak{k}$. If $\Gamma(t) = exp(tX)(o)$ is a geodesic for some $X \in \mathfrak{m}$, then $\Gamma(t)$ is also a geodesic for the canonical connection $\tilde{\nabla}$. For this reason, such geodesics are called {\em canonically homogeneous}. In general, a homogeneous geodesic is not canonically homogeneous. Note that any reductive homogeneous Lorentzian manifold $(M,g)$ admits at least one homogeneous geodesic through a point [P]. The question whether such a space always admits a {\em null} homogeneous geodesic through a point, had a negative answer in [CM2]. Further examples will be given in Section 3.

We now recall how the geometric problem of finding homogeneous geodesics of a reductive homogeneous space, reduces to the algebraic problem of determining its geodesic vectors. 
Let $(M=K/H,g)$ be a reductive homogeneous Lorentzian manifold and $\mathfrak{k}$ = $\mathfrak{m} \oplus \mathfrak{h}$ the corresponging reductive split of the Lie algebra $\mathfrak{k}$. The canonical projection $ p : K
\rightarrow  K/H $ induces an isomorphism between the subspace $\mathfrak{m}$
 and the tangent space $T_{o}(M)$. In particular, the Lorentzian metric $g_{o}$ on $T_{o}(M)$ induces a Lorentzian metric
$<,>$ on $\mathfrak{m}$, which is Ad($H$)-invariant. The following characterization is a crucial step for determining homogeneous geodesics of a reductive homogeneous pseudo-Riemannian manifold:

\begin{proposition} {\bf ([P],[FMP],[DK]).} Consider a geodesic $\Gamma (t)$ of $M=K/H$, with $\Gamma (0)=o$ and $\Gamma ' (0)= X_{\mathfrak{m}} \in \mathfrak{m} \, (\equiv T_o (K/H))$. $\Gamma (t)$ is homogeneous if and only if there exists 
$X_{\mathfrak{h}} \in \mathfrak{h}$ such that $X=X_{\mathfrak{m}} + X_{\mathfrak{h}} \in \mathfrak{k}$ satisfies
\begin {equation}
<[X,Y]_{\mathfrak{m}},  X_{\mathfrak{m}} > = k <X_{\mathfrak{m}},  Y >,
\end{equation}
for all $Y\in \mathfrak{m}$ and some $k \in \Real$ depending on $X_{\mathfrak{m}}$. 
\end{proposition}
Proposition 2.3, whose proof can be found in [DK], is the Lorentzian analogue of Proposition 2.1 of [KV], characterizing homogeneous geodesics of a Riemannian homogeneous space.

A vector $X \in \mathfrak{k}$ satisfying (2.3) is called a {\em geodesic vector}.
When $X_{\mathfrak{m}}$ is either spacelike ($<X_{\mathfrak{m}},X_{\mathfrak{m}}> >0$) or timelike ($<X_{\mathfrak{m}},X_{\mathfrak{m}}> <0$), applying (2.3) with $Y=X_{\mathfrak{m}}$ we get $k=0$, while for a null vector $X_{\mathfrak{m}}$, $k$ may be any real constant. Note also that if $\mathfrak{h}=0$, then $\mathfrak{k}=\mathfrak{m}$ and (2.3) simplifies as follows:
\begin {equation}
<[X,Y], X > = k <X,  Y>,
\end{equation}
for all $Y \in \mathfrak{k}$.

A finite family $\{\Gamma_1,\Gamma_2,...,\Gamma_k\}$ of
homogeneous geodesics through ${o}\in M =K/H$ is said to be {\em linearly independent}
if the corresponding initial tangent vectors at $o$ are linearly independent. 
The following result is obvious.

\begin{proposition} A finite family
$\{\Gamma_1,\Gamma_2,...,\Gamma_k\}$ of homogeneous geodesics
through ${p_o}\in M =K/H$ is linearly
independent if the $\mathfrak{m}$-components of the corresponding geodesic
vectors are linearly independent.
\end{proposition}

\bigskip \noindent
\section{\red Three-dimensional unimodular Lorentzian Lie groups and their \\ homogeneous geodesics}
\setcounter{equation}{0}

S. Rahmani [R] classified three-dimensional unimodular Lie groups equipped with a left-invariant Lorentzian metric, obtaining a result corresponding to the one found by Milnor [Mi] in the Riemannian case. Earlier, Cordero and Parker [CoPa2] already studied three-dimensional Lie groups equipped with left-invariant Lorentzian metrics, determining their curvature tensors and investigating the symmetry groups of the sectional curvature in the different cases. In particular, they wrote down the possible forms of a non-unimodular Lie algebra.  Taking into account these results and Theorem 1.1, we have the following:

\begin{theorem}{\bf [C]}
Let $(M,g)$ be a three-dimensional connected, simply connected, complete homogeneous Lorentzian manifold. If $(M,g)$ is not symmetric, then $M=G$ is a three-dimensional Lie group and $g$ is left-invariant. Precisely, one of the following cases occurs:
\begin{itemize}
\item
If $G$ is unimodular, then there exists a pseudo-orthonormal frame field $\{ e_1,e_2,e_3\}$, with $e_3$ timelike, such that the Lie algebra of $G$ is one of the following:

\medskip
 a) 
 \begin{eqnarray}
& &\left[e_1,e_2 \right]=\alpha e_1-\beta e_3, \nonumber \\ 
(\mathfrak{g} _1 ) : & &\left[ e_1,e_3\right]=-\alpha e_1-\beta e_2, \\
& & \left[e_2,e_3\right]=\beta e_1 +\alpha e_2 +\alpha  e_3 \qquad \alpha \neq 0.  \nonumber
\end{eqnarray}   
In this case, $G= O(1,2)$ or $SL(2,\Real)$ if $\beta \neq 0$, while $G=E(1,1)$ if $\beta=0$.

\medskip
 b) 
\begin{eqnarray}
& &\left[e_1,e_2 \right]=\gamma e_2-\beta e_3, \nonumber \\ 
(\mathfrak{g} _2): & &\left[ e_1,e_3\right]=-\beta e_2+\gamma e_3, \qquad \gamma \neq 0, \\
& & \left[e_2,e_3\right]=\alpha e_1 .  \nonumber
\end{eqnarray}   
In this case, $G= O(1,2)$ or $SL(2,\Real)$ if $\alpha \neq 0$, while $G=E(1,1)$ if $\alpha=0$.

\medskip
c) 
\begin{eqnarray}
& &\left[e_1,e_2 \right]=-\gamma e_3, \nonumber \\ 
(\mathfrak{g} _3): & &\left[ e_1,e_3\right]=-\beta e_2, \\
& & \left[e_2,e_3\right]=\alpha e_1 .  \nonumber
\end{eqnarray}   
{\em The following Table I lists all the Lie groups $G$ which admit a Lie algebra {\bf $g_3$}, taking into account the different possibilities for $\alpha$, $\beta$ and $\gamma$:}
\begin{center}
\begin{tabular}{|c|c|c|c|}
\hline $G$ &  $\alpha$ & $\beta$ & $\gamma$ \\
\hline  $O(1,2)$ or $SL(2,\Real)$ & $+$ & $+$ & $+$\\
\hline  $O(1,2)$ or $SL(2,\Real)$ & $+$ & $-$ & $-$\\
\hline  $SO(3)$ or $SU(2)$ & $+$ & $+$ & $-$\\
\hline  $E(2)$ & $+$ & $+$ & $0$\\
\hline  $E(2)$ & $+$ & $0$ & $-$\\
\hline  $E(1,1)$ & $+$ & $-$ & $0$\\
\hline  $E(1,1)$ & $+$ & $0$ & $+$\\
\hline  $H_3$ & $+$ & $0$ & $0$\\
\hline  $H_3$ & $0$ & $0$ & $-$\\
\hline  $\Real \oplus \Real \oplus \Real$ & $0$ & $0$ & $0$\\
\hline 
\end{tabular} \nopagebreak \\ \nopagebreak Table I $\vphantom{\displaystyle\frac{a}{2}}$
\end{center}

\medskip
d) 
 \begin{eqnarray}
& &\left[e_1,e_2 \right]=- e_2 + (2 \varepsilon - \beta) e_3, \qquad \varepsilon = \pm 1, \nonumber \\ 
(\mathfrak{g} _4): & &\left[ e_1,e_3\right]=-\beta e_2 + e_3,  \\
& & \left[e_2,e_3\right]=\alpha e_1 .  \nonumber
\end{eqnarray}   
{\em The following Table II describes all Lie groups $G$ admitting a Lie algebra {\bf $g_4$}:}
\begin{gather*}
\begin{array}{cc}
\begin{tabular}{|c|c|c|}
\hline $G$\quad  {\rm ($\varepsilon =1$)} &  $\alpha$ & $\beta$  \\ 
\hline  $O(1,2)$ or $SL(2,\Real)$ & $\neq 0$ & $\neq 1$ \\
\hline  $E(1,1)$ & $0$ & $\neq 1$ \\
\hline  $E(1,1)$ & $<0$ & $1$ \\
\hline  $E(2)$ & $>0$ & $1$ \\
\hline  $H_3$ & $0$ & $1$ \\
\hline
\end{tabular} 
& \qquad
\begin{tabular}{|c|c|c|}
\hline $G$ \quad {\rm ($\varepsilon =-1$)} &  $\alpha$ & $\beta$  \\ 
\hline  $O(1,2)$ or $SL(2,\Real)$ & $\neq 0$ & $\neq -1$ \\
\hline  $E(1,1)$ & $0$ & $\neq -1$ \\
\hline  $E(1,1)$ & $>0$ & $-1$ \\
\hline  $E(2)$ & $<0$ & $-1$ \\
\hline  $H_3$ & $0$ & $-1$ \\
\hline
\end{tabular} 
\end{array}
 \\ 
 \text{Table II} \vphantom{\displaystyle\frac{a}{2}}
\end{gather*}

\item
If $G$ is non-unimodular, then there exists a pseudo-orthonormal frame field $\{ e_1,e_2,e_3\}$, with $e_3$ timelike, such that the Lie algebra of $G$ is one of the following: 

\medskip
e) 
\begin{eqnarray}
& &\left[e_1,e_2 \right]=0, \nonumber \\ 
(\mathfrak{g} _5): & &\left[ e_1,e_3\right]=\alpha e_1+\beta e_2, \\
& & \left[e_2,e_3\right]=\gamma e_1 +\delta e_2, \qquad \alpha +\delta \neq 0, \, 
\alpha \gamma +\beta \delta =0.  \nonumber
\end{eqnarray}   
\medskip
f) 

 \begin{eqnarray}
& &\left[e_1,e_2 \right]=\alpha e_2 +\beta e_3, \nonumber \\ 
(\mathfrak{g} _6): & &\left[ e_1,e_3\right]=\gamma e_2+\delta e_3,  \\
& & \left[e_2,e_3\right]= 0, \qquad \qquad \qquad \alpha +\delta \neq 0, \, \alpha \gamma -\beta \delta =0.  \nonumber
\end{eqnarray}   
\medskip
 g) 
\begin{eqnarray}
& &\left[e_1,e_2 \right]=-\alpha e_1-\beta e_2 -\beta e_3, \nonumber \\ 
(\mathfrak{g} _7): & &\left[ e_1,e_3\right]=\alpha e_1+\beta e_2 +\beta e_3, \\
& & \left[e_2,e_3\right]=\gamma e_1 +\delta e_2 +\delta  e_3 ,  \qquad \alpha +\delta \neq 0, \, \alpha \gamma =0. \nonumber
\end{eqnarray}   
\end{itemize}

\end{theorem}
  
Following [CoPa2], cases $(\mathfrak{g} _5)-(\mathfrak{g} _7)$ are the possible forms of the non-unimodular Lie algebra of a  three-dimensional Lorentzian Lie group, rewritten here for a Lorentzian metric of signature $(+,+,-)$ and a pseudo-orthonormal frame field $\{ e_1,e_2,e_3\}$ with $e_3$ timelike. The determinant $D=\frac{4(\alpha \delta -\beta \gamma)}{(\alpha +\delta )^2}$ provides a complete isomorphism invariant for Lie algebras $(\mathfrak{g} _5)-(\mathfrak{g} _7)$.

\medskip
Theorems 1.1 and 3.1 above have been used in [C] also to obtain the full classification of three-dimensional Lorentzian symmetric spaces. The results are summarized in the following

\begin{theorem}{\bf [C]}
A connected, simply connected three-dimensional Lorentzian symmetric space $(M,g)$ is either 

\medskip
i) a Lorentzian space form $S ^3 _1$, $\Real ^3 _1$ or $\Hyp ^3 _1$, or

\medskip
ii) a direct product $\Real \times S ^2 _1$, $\Real \times \Hyp ^2 _1$, $S^2 \times 
\Real $ or $\Hyp ^2 \times \Real$, or

\medskip
iii) a symmetric Lorentzian manifold having a {\em parallel null vector field}. It admits local coordinates $(t,x,y)$ such that, with respect to the local frame field $\{ (\frac{\partial}{\partial t}), (\frac{\partial}{\partial x}), (\frac{\partial}{\partial y}) \}$, the Lorentzian metric $g$ and the Ricci operator are given by 
\begin{equation}
g=
\left(
\begin{array}{ccc}
0 &  0 & 1   \\
0 & \varepsilon & 0 \\
1 & 0 & f
\end{array}
\right),  \qquad
 Q=
\left(
\begin{array}{ccc}
0 &  0 & -\frac{1}{\varepsilon} \alpha   \vphantom{\displaystyle\frac{a}{2}} \\
0 & 0 & 0 \vphantom{\displaystyle\frac{a}{2}} \\
0 & 0 & 0
\end{array}
\right),
\end{equation}   
where $\varepsilon = \pm 1$, $u = (\frac{\partial}{\partial t})$ and $
f(x,y)=x^2 \alpha +x \beta (y) + \xi (y)$, for any constant $\alpha \in \Real$ and any functions $\beta,\xi$ {\rm [ChGV, Theorem 6]}. 
\end{theorem}

In particular, from the proof of Theorem 3.2, the following result holds for a non-unimodular Lie group:

\begin{proposition}
A three-dimensional non-unimodular Lie group $G$, equipped with a left-inva- riant Lorentzian metric $g$, is symmetric if and only if its Lie algebra is one of the following:
\begin{itemize}
\item
$\mathfrak{g_5}$ with either $\alpha =\beta = \gamma =0 \neq \delta$, $\beta = \gamma = \delta =0 \neq \alpha$  or $\beta + \gamma =0 \neq \alpha =\delta$. 
\item
$\mathfrak{g_6}$ with either $\alpha =\beta = \gamma =0 \neq \delta$, $\beta = \gamma = \delta =0 \neq \alpha$, $\beta - \gamma =0 \neq \alpha =\delta$ or $ \beta - \varepsilon \alpha =0=\gamma - \varepsilon \delta$, with $\varepsilon = \pm 1$. 
\item
$\mathfrak{g_7}$ with either  $\alpha = \gamma =0\neq \delta$,  $\gamma =\delta = 0 \neq \alpha$ or  $\alpha -\delta = \gamma =0$. 
\end{itemize}

\end{proposition}

Next, let $(G,g)$ be a three-dimensional non-unimodular Lie group equipped with a left-invariant Lorentzian metric $g$; let $\mathfrak{g}$ its Lie algebra.  As group space, G admits a left-invariant connection $\tilde{\nabla}$, such that $\tilde{\nabla} _X Y=0$ for all $X,Y \in \mathfrak{g}$ [KoNo]. Let $\tilde{R}$ and $\tilde{T}$ denote the curvature tensor and the torsion tensor associated to $\tilde{\nabla}$, respectively. Then, $\tilde{R}(X,Y)=0$ and $\tilde{T} (X,Y) = -[X,Y]$ for all vector fields $X,Y \in \mathfrak{g}$. By Theorem 2.2, the tensor $T = \nabla - \tilde{\nabla}$, where $\nabla$ is the Levi-Civita connection, is a Lorentzian homogeneous structure.

For all possible forms $\mathfrak{g} _i$, $i=5,6,7$, we shall compute the Lie algebra $\mathfrak{k}$ of the connected component $K$ of the full isometry group $I(G)$, find a reductive decomposition and determine geodesic vectors.       

In order to determine $\mathfrak{h}$ so that $\mathfrak{k} = \mathfrak{g} \oplus \mathfrak{h}$, we first note that $\mathfrak{h}$ is isomorphic to the set of all self-adjoint endomorphisms $A$ of $\mathfrak{g}$, such that $A(R)=A(\nabla R)=...=A(\nabla ^n R)=0$ for all $n$. Because $G$ is three-dimensional, we  equivalently have 

\begin{displaymath}
\mathfrak{h} =\{ A \in gl (\mathfrak{g}) / A(g)=A(\varrho)=A(\nabla \varrho)=...=A(\nabla ^n \varrho)...=0 \}.
\end{displaymath}  
Moreover, put $$\mathfrak{h}_k =\{ A \in gl (\mathfrak{g}) / A(g)=A(\varrho)=A(\nabla \varrho)=...=A(\nabla ^k \varrho)=0 \},$$ 
for any $k=0,1,..$, and $$\mathfrak{l} =\{ A \in gl (\mathfrak{g}) / A(g)=A(\tilde{R})=A(\tilde{T})=0 \}.$$ 
Then, we clearly have
\begin{displaymath}
\mathfrak{l} \subset \mathfrak{h} \subset \mathfrak{h}_k ,
\end{displaymath}   
for all $k$ and so, if $\mathfrak{l} = \mathfrak{h}_k$ for some $k$, then $\mathfrak{h} = \mathfrak{l}$. We refer to [K1, p.54] for more details.  

Let $A$ be an endomorphism of $\mathfrak{g}$ and $\{e_1,e_2,e_3\}$ a pseudo-orthonormal basis of $\mathfrak{g}$, with $e_3$ timelike. Condition $A(g)=0$ means that $A$ is self-adjoint. Hence, there exist some real constants $a,b,c$, such that
 \begin{equation}\label{metr}
A e_1 = b e_2 + c e_3, \quad A e_2 = -b e_1 + a e_3, \quad A e_3 = c e_1 + a e_2. 
\end{equation}   

Since $\tilde{R}=0$ $\tilde{T} (X,Y) = -[X,Y]$, condition $A(\tilde{R})=0$ is trivially satisfied, while $A(\tilde{T})=0$ is equivalent to
\begin{equation}\label{parlie}
A \left[  e_i,e_j\right]= \left[ A e_i,e_j \right] + \left[  e_i, A e_j \right] \quad {\rm for \; all} \; i,j.
\end{equation}   

Now, we shall inspect our three basic algebras case by case.

\medskip
$(\mathfrak{g} _5 )$:  Consider an endomporphism $A$ of $\mathfrak{g}_5$ satisfying (\ref{metr}). Starting from (3.5), we can write down (\ref{parlie}) and we get
\begin{equation}\label{sys1}
\left\{
\begin{array}{l}
\alpha a -\gamma c=0, \\
\beta a -\delta c= 0, \\
(\beta +\gamma)b=0, \\
(\alpha -\delta)b=0, \\
\beta a +\alpha c=0, \\
\delta a +\gamma c=0.
\end{array}
\right.
\end{equation}
Taking into account $\alpha +\delta \neq 0$ and $\alpha \gamma +\beta \delta =0$, we have two possible cases:

\medskip
(I): $(\gamma, \delta) \neq (-\beta, \alpha)$. In this case, $a=b=c=0$ is the only solution of (\ref{sys1}), that is, $\mathfrak{l}=0$.

\medskip
(II): $(\gamma, \delta) = (-\beta, \alpha)$. Then, by (\ref{sys1}), we only get $a=c=0$. So, $\mathfrak{l}={\rm Span}(A)$, where $A e_1 =e_2$, $A e_2 = -e_1$ and $A e_3 =0$.

\medskip
$(\mathfrak{g} _6 )$:  Let $A$ be an endomporphism of $\mathfrak{g}_6$ satisfying (\ref{metr}). We use (3.6) to compute (\ref{parlie}) and we obtain
\begin{equation}\label{sys2}
\left\{
\begin{array}{l}
 (\alpha-\delta)a =0, \\
(\beta -\gamma )a=0, \\
\alpha b -\beta c =0, \\
\gamma b -\delta c =0, \\
\gamma b + \alpha c =0, \\
\delta b + \beta c =0.
\end{array}
\right.
\end{equation}
Since $\alpha +\delta \neq 0$ and $\alpha \gamma -\beta \delta =0$, we have that if $(\gamma, \delta) \neq (\beta, \alpha)$, then $a=b=c=0$ is the only solution of (\ref{sys2}). Therefore, in this case $\mathfrak{l}=0$. In the remaining case, $(\gamma, \delta) = (\beta, \alpha)$ and so, by Proposition 3.3, $(G,g)$ is symmetric. 

\medskip
$(\mathfrak{g} _7 )$:  Let $A$ be an endomporphism of $\mathfrak{g}_7$, satisfying (\ref{metr}). By (3.7), we see that $A$ satisfies (\ref{parlie}) if and only if 
\begin{equation}\label{sys3}
\left\{
\begin{array}{l}
\beta (b-c)=\alpha a - \gamma c, \\
\beta a + \alpha b =\delta c -\beta a , \\
\beta a + \alpha c =\delta c -\beta a, \\
\beta(b -c) =\alpha a -\gamma b, \\
\beta a + \alpha b =\delta b -\beta a , \\
\beta a + \alpha c =\delta b -\beta a, \\
\alpha (c-b) =\delta (c-b), \\
\beta (c- b)= \delta a +\gamma b, \\
\beta (c- b)= \delta a +\gamma c. \\
\end{array}
\right.
\end{equation}
It is easy to show that, whenever $\alpha \neq \delta$, we have 
$a=b=c=0$ as the only solution of (\ref{sys3}) and so, $\mathfrak{l}=0$. On the other hand, if $\alpha =\delta$ then, by (3.7), $\alpha \neq 0$ and $\gamma =0$. In this case, Proposition 3.3 implies that $(G,g)$ is symmetric.

\medskip
Next, routine but very long calculations can show that, for any $\mathfrak{g} =\mathfrak{g} _i$, $i=5,6,7$, there exists a $k \leq 2$ such that $\mathfrak{l} = \mathfrak{h}_k$. In particular, we have the following

\begin{proposition} 
{\em Let $(G,g)$ be a non-unimodular Lie group, equipped with a {\em non-symmetric} left-invariant Lorentzian metric, and $\mathfrak{g}$ its Lie algebra.
\begin{itemize}
\item When $\mathfrak{g}=\mathfrak{g}_5$, we have $\mathfrak{h}_0 =\mathfrak{l}$, except in the following cases: 

\medskip
a) if $(\gamma,\delta) \neq (-\beta, \gamma )$ and $\beta \delta \neq 0$, then $\mathfrak{h}_1 =\mathfrak{l}=0$. 

\medskip
b) if either $\alpha =\beta =0$ and $\gamma \neq 0 \neq  \delta$, or $\gamma =\delta =0$ and $\alpha \neq 0 \neq \beta$, then $\mathfrak{h}_2 =\mathfrak{l}=0$. 

\item When $\mathfrak{g}=\mathfrak{g}_6$, we have $\mathfrak{h}_0 =\mathfrak{l}=0$, unless  $(\gamma,\delta) \neq (\beta, \gamma )$ and $\beta(\beta ^2 -\alpha ^2)\neq 0$. In the last case, $\mathfrak{h}_2 =\mathfrak{l}=0$. 

\item When $\mathfrak{g}=\mathfrak{g}_7$, we have $\mathfrak{h}_0 =\mathfrak{l}=0$, except in the following cases:

\medskip
a) if $\gamma=0$ and $\alpha \delta (\alpha ^2 -\delta ^2)\neq 0$, then $\mathfrak{h}_1 =\mathfrak{l}=0$. 

\medskip
b) if $\alpha =\beta =0\neq \gamma$, then $\mathfrak{h}_2 =\mathfrak{l}=0$. 
\end{itemize} }
\end{proposition}

\noindent
By Proposition 3.4, we have $\mathfrak{h}=\mathfrak{l}$ for all non-symmetric non-unimodular Lorentzian Lie groups. Therefore, for each of them, $\mathfrak{k}=\mathfrak{g} \oplus \mathfrak{l}$ is the Lie algebra of the connected component $K$ of the full isometry group of $G$. So, we can now determine the homogeneous geodesics in the different cases.

\medskip
$(\mathfrak{g} _5 )$:  If $(\gamma,\delta) = (-\beta, \alpha )$ then, by Proposition 3.3, $(G,g)$ is symmetric. In particular, it is naturally reductive. Hence, in the sequel we shall assume $(\gamma,\delta) \neq (-\beta, \gamma )$. In this case, $\mathfrak{h}=\mathfrak{l}=0$ and so, $\mathfrak{k}=\mathfrak{g}_5$ and $K=G$. 
A vector $X \in \mathfrak{g}_5$ is geodesic if and only if (2.4) holds for all $Y \in \mathfrak{g}_5$ and some constant $k$. We write down (2.4) for $Y=e_i$, $i=1,2,3$. Using (3.5) and taking into account that $X = x_1 e_1 +x_2 e_2 +x_3 e_3$, we find that (2.4) is equivalent to the following system: 
\begin{equation}
\left\{
\begin{array}{l}
-x_3(\alpha x_1+ \beta x_2) =k x_1, \vphantom{\displaystyle\frac{a}{2}} \\
- x_3 (\gamma x_1 + \delta x_2)= k x_2, \vphantom{\displaystyle\frac{a}{2}} \\
 \alpha x_1 ^2 + (\beta +\gamma)x_1 x_2 +\delta x_2 ^2=-k x_3. \vphantom{\displaystyle\frac{a}{2}} \\
\end{array}
\right.
\end{equation}
In determining the solutions of (3.14) we must also take into account that, by (3.5), $\alpha +\delta \neq 0$ and $\alpha \gamma +\beta \delta =0$. It is easy to prove that when $k =0$, the solutions of (3.14) are either $x_3=\alpha x_1 ^2+ (\beta +\gamma)x_1 x_2+\delta x_2 ^2=0$,  $x_1=x_2=0$, $\alpha x_1 +\beta x_2 =0$ (but only if $\gamma=\delta=0$), or $\gamma x_1 +\delta x_2 =0$ (but only if $\alpha=\beta=0$). If $k \neq 0$, only null-like solutions can occur, so we must add we must add to (3.14) the condition 
\begin{equation}\label{null}
x_1 ^2 + x_2 ^2 - x_3 ^2 =0. 
\end{equation}
Standard calculations show that, when $k \neq 0$, (3.14) and (3.15) are both satisfied if and only if $\gamma x_1 ^2+ (\delta -\alpha)x_1 x_2-\beta x_2 ^2=0$. Summarizing,  geodesic vectors have one of the following  forms:
\begin{equation}
\left\{
\begin{array}{l}
 X= x_1 e_1 +x_2 e_2, \quad {\rm with} \; \alpha x_1 ^2+ (\beta +\gamma)x_1 x_2+\delta x_2 ^2=0, \vphantom{\displaystyle\frac{a}{2}} \\
 X= x_3 e_3 , \vphantom{\displaystyle\frac{a}{2}} \\
 X= x_1 ( \delta e_1 -\gamma e_2) +x_3 e_3 \quad {\rm but } \;
\alpha =\beta=0,\vphantom{\displaystyle\frac{a}{2}}  \\
 X= x_2 ( -\beta e_1 +\alpha e_2) +x_3 e_3 \quad {\rm but } \;
\gamma =\delta=0, \vphantom{\displaystyle\frac{a}{2}} \\
 X= x_1 e_1 +x_2 e_2 \pm \sqrt{x_1^2 +x_2 ^2} \, e_3 \quad {\rm  \, with} \; \gamma x_1 ^2+ (\delta -\alpha)x_1 x_2-\beta x_2 ^2=0. \vphantom{\displaystyle\frac{a}{2}}
\end{array}
\right. 
\end{equation}
Note that, by (3.16), whenever $(\alpha,\beta)\neq (0,0) \neq (\gamma,\delta)$ and $(\alpha -\delta)^2 +4\beta \gamma <0$, there are not null-like geodesic vectors. So, in this case $G$ does not admit null homogeneous geodesics. If in addition $(\beta +\gamma)^2 -4\alpha \delta <0$, then the only geodesic vectors are the ones parallel to $e_3$ (timelike).

\medskip
$(\mathfrak{g} _6 )$:  We assume $(\gamma,\delta) \neq (\beta, \gamma )$, since in the remaining case, $G$ is symmetric. When $(\gamma,\delta) \neq (\beta, \gamma )$, we have $\mathfrak{h}=\mathfrak{l}=0$ and so, $\mathfrak{k}=\mathfrak{g}_6$ and $K=G$. 
This case is quite similar to the corresponding one for $\mathfrak{g}_5$. In fact, a vector $X = x_1 e_1 + x_2 e_2 + x_3 e_3 \in \mathfrak{k}$ is geodesic if and only if satisfies (2.4), that is, using (3.6), 
\begin{equation}
\left\{
\begin{array}{l}
-\alpha x_2 ^2 + (\beta -\gamma)x_2 x_3 +\delta x_3 ^2=k x_1, \vphantom{\displaystyle\frac{a}{2}} \\
x_1(\alpha x_2 -\beta x_3) =k x_2,  \vphantom{\displaystyle\frac{a}{2}} \\
x_1 (\gamma x_2- \delta x_3)= -k x_3. \vphantom{\displaystyle\frac{a}{2}} \\
\end{array}
\right.
\end{equation}
By (3.5) we have $\alpha +\delta \neq 0$ and $\alpha \gamma -\beta \delta =0$. For $k =0$, the solutions of (3.19) are either $x_1=\alpha x_2 ^2+ (\gamma -\beta)x_2 x_3-\delta x_3 ^2=0$,  $x_2=x_3=0$, $\alpha x_2-\beta x_3 =0$ (but only if $\gamma=\delta=0$), or $\gamma x_2 -\delta x_3 =0$ (but only if $\alpha=\beta=0$). Indeed, $\alpha x_2-\beta x_3 =0$ and $\gamma x_2 -\delta x_3 =0$ also occur as solutions of (3.19) when $(\beta,\delta)=\pm (\alpha,\gamma)$, but in this case we get a symmetric space [C] and so, we discard it.

If $k \neq 0$, the null-like solutions of (3.19) are all and the ones triples $(x_1,x_2,x_3)$ satisfying $\gamma x_2 ^2+ (\alpha -\delta)x_2 x_3-\beta x_3 ^2=0$ (and the null-like condition(\ref{null})). Therefore, geodesic vectors have one of the following  forms:
\begin{equation}
\left\{
\begin{array}{l}
 X= x_2 e_2 +x_3 e_3, \quad {\rm with} \; \alpha x_1 ^2+ (\gamma-\beta)x_2 x_3-\delta x_3 ^2=0, \vphantom{\displaystyle\frac{a}{2}}\\
 X= x_1 e_1 , \vphantom{\displaystyle\frac{a}{2}} \\
 X= x_1 e_1 + x_2 ( \delta e_2 +\gamma e_3) \quad {\rm but } \;
\alpha =\beta=0, \vphantom{\displaystyle\frac{a}{2}} \\
 X= x_1 e_1 + x_3 ( \beta e_2 +\alpha e_3)  \quad {\rm but } \;
\gamma =\delta=0, \vphantom{\displaystyle\frac{a}{2}} \\
 X= x_1 e_1 +x_2 e_2 \pm \sqrt{x_1^2 +x_2 ^2} \, e_3 \quad {\rm  \, with} \; \gamma x_2 ^2+ (\alpha-\delta)x_2 x_3-\beta x_3 ^2=0. 
\end{array}
\right. 
\end{equation}
By (3.20), whenever $(\beta -\gamma)^2 +4\alpha \delta <0$, $(\alpha,\beta)\neq (0,0) \neq (\gamma,\delta)$ and $(\alpha -\delta)^2 +4\beta \gamma <0$, the only geodesic vectors are the ones parallel to $e_1$ (spacelike).

\medskip
$(\mathfrak{g} _7 )$: By (3.7), $\alpha \gamma =0$, that is, either $\alpha =0$ or $\gamma =0$. It is enough to consider the case $\alpha \neq \delta$, since in the remaining case  we have a symmetric space. In fact, if $\alpha =\delta$, then (3.7) implies $2 \alpha = \alpha +\delta  \neq 0$ and so, $\gamma =0$. Then, $\alpha -\delta =\gamma =0$ and, by Proposition 3.3, $(G,g)$ is symmetric. 

So, in the sequel we shall assume $\alpha \neq \delta$. Then, $\mathfrak{h}=0$, that is, $\mathfrak{k}=\mathfrak{g}_7$ and $K=G$. Applying (2.4) and taking into account (3.7), we find that a vector $X \in \mathfrak{g}_7$ is geodesic if and only if 
\begin{equation}
\left\{
\begin{array}{l}
(x_2-x_3)[\alpha x_1 +\beta (x_2-x_3)]=k x_1, \vphantom{\displaystyle\frac{a}{2}} \\
(\alpha x_1 + \gamma x_3) x_1 + (\beta x_1 + \delta x_3) (x_2- x_3) = -k x_2, \vphantom{\displaystyle\frac{a}{2}} \\
(\alpha x_1 + \gamma x_2) x_1 + (\beta x_1 + \delta x_2) (x_2- x_3) = -k x_3. \vphantom{\displaystyle\frac{a}{2}}
\end{array}
\right.
\end{equation}
We exclude the case $\alpha=\gamma=0$, since by Proposition 3.3 it corresponds to a symmetric space. In order to solve system (3.21), we shall treat separately the cases $\alpha =0 \neq \gamma$ and $\gamma =0 \neq \alpha$.

\medskip A): $\alpha =0 \neq \gamma$. Then, (3.21) becomes
\begin{equation}
\left\{
\begin{array}{l}
\beta (x_2-x_3)^2=k x_1, \vphantom{\displaystyle\frac{a}{2}} \\
\gamma x_1 x_3 + (\beta x_1 + \delta x_3) (x_2- x_3) =- k x_2, \vphantom{\displaystyle\frac{a}{2}} \\
\gamma x_1 x_2 + (\beta x_1 + \delta x_2) (x_2- x_3) = -k x_3. \vphantom{\displaystyle\frac{a}{2}}
\end{array}
\right.
\end{equation}
By standard calculations we get that when $k=0$, (3.22) holds if and only if either $x_2=x_3=0$, $x_1=x_2-x_3=0$ or $\gamma x_1+\delta (x_2-x_3)=0$ (the latter only if $\beta=0$). When $k\neq 0$, we must consider only null-like solutions. We treat separately the cases $x_1 =0$ and $x_1 \neq 0$. If $x_1 =0$, we easily conclude that (3.22) only admits the solution $x_1 =x_2+x_3=0$, and only if $\beta =0$. In the remaining case $x_1 \neq 0$, we can use the first equation of (3.22) to write $k$ in function of $x_1$, $x_2$ and $x_3$. Then, rather long calculations lead to show that the only null-like solutions of (3.22) (with $x_1 \neq 0$) are of the form 
$$x_1 = \frac{\beta +\gamma}{\delta} x_2 -\frac{\beta -\gamma}{\delta}x_3,$$
where $x_2$, $x_3$ must satisfy the second order homogeneous equation
\begin{equation}\label{cond}
[(\beta +\gamma)^2 +\delta ^2]x_2 ^2 -2 (\beta ^2 -\gamma ^2) x_2 x_3 +[(\beta -\gamma)^2 -\delta ^2]x_3 ^2 =0.
\end{equation}
Equation (\ref{cond}) only admits real solutions when $\delta ^2 (4 \beta \gamma +\delta ^2) \geq 0$. Since $\alpha =0$, we have $D=\frac{-4\beta\gamma}{\delta ^2}$. Therefore, $\delta ^2 (4 \beta \gamma +\delta ^2) =\delta ^4 (1-D)$ and so, such solutions exist if and only if $D \leq 1$. Hence, all geodesic vectors are of one of the following forms:
\begin{equation}
\left\{
\begin{array}{l}
X=x_1 e_, \vphantom{\displaystyle\frac{a}{2}} \\
 X= x_2 (e_2 + e_3), \vphantom{\displaystyle\frac{a}{2}}\\
 X=  -\frac{\delta}{\gamma} (x_2-x_3) e_1 +x_2 e_2 + x_3 e_3 \quad {\rm but } \;
\beta=0, \vphantom{\displaystyle\frac{a}{2}} \\
X= \left(\frac{\beta + \gamma}{\delta} x_2 -  \frac{\beta -\gamma}{\delta} x_3 \right) e_1 +x_2 e_2 + x_3 e_3 \; {\rm but \; (3.23) \; holds } . \vphantom{\displaystyle\frac{a}{2}} 
\end{array}
\right. 
\end{equation}

\medskip B): $\alpha \neq 0 =\gamma$. In this case, (3.21) becomes
\begin{equation}
\left\{
\begin{array}{l}
(x_2-x_3)[\alpha x_1 +\beta (x_2-x_3)]=k x_1, \vphantom{\displaystyle\frac{a}{2}} \\
\alpha x_1 ^2 + (\beta x_1 + \delta x_3) (x_2- x_3) = -k x_2, \vphantom{\displaystyle\frac{a}{2}} \\
\alpha x_1 ^2 + (\beta x_1 + \delta x_2) (x_2- x_3) = -k x_3. \vphantom{\displaystyle\frac{a}{2}}
\end{array}
\right.
\end{equation}
Suppose first that $k=0$. Then, it is easy to conclude that, apart from cases corresponding to symmetric spaces listed in Proposition 3.3, (3.25) holds if and only if  $x_1=x_2-x_3=0$. For $k\neq0$, if $x_1 =0$, then (3.25) admit the solution $x_1=x_2+x_3=0$, but only if $\beta =0$. In the remaining case $x_1 \neq 0$, after some calculations we obtain that the solutions are all and the ones of the form  
\begin{equation}
\left\{
\begin{array}{l}
x_1 = \frac{2\beta (\alpha -\delta)}{\beta ^2 +(\alpha -\delta )^2} x_3, \vphantom{\displaystyle\frac{a}{2}} \\
x_2 = \frac{\beta ^2 -(\alpha -\delta )^2}{\beta ^2 +(\alpha -\delta )^2} x_3,\vphantom{\displaystyle\frac{a}{2}} \\
\end{array}
\right.
\end{equation}
Therefore, all geodesic vectors are of one of the following forms:
\begin{equation}
\left\{
\begin{array}{l}
 X= x_2 (e_2 + e_3), \vphantom{\displaystyle\frac{a}{2}}\\
 X=  x_2( e_2 - e_3) \quad {\rm but } \; \beta=0, \vphantom{\displaystyle\frac{a}{2}} \\
X= x_3 \left\{ 2\beta (\alpha -\delta ) e_1 + [\beta ^2 -(\alpha -\delta )^2] e_2 +
[\beta ^2 +(\alpha -\delta )^2] e_3 \right\}. \vphantom{\displaystyle\frac{a}{2}} 
\end{array}
\right. 
\end{equation}

The calculations above are resumed in the following

\begin{proposition}
Let $(M,g)$ be a three-dimensional Lorentzian homogeneous space, not symmetric, isometric to a non-unimodular Lie group $G$. Let $\mathfrak{g} _i$, $i=5,6,7$, denote the Lie algebra of $G$. Then, the set of geodesic vectors of $M$ through any point is described in the following {\rm Table III}.  
\end{proposition}

\begin{center}
\begin{tabular}{|c|l|}
\hline Table III &  \\ 
\hline {\rm Lie algebra} &  {\rm Geodesic vectors} \\ & \\

\hline  $\mathfrak{g} _5 $ with $(\gamma,\delta) \neq (-\beta,\alpha)$ $\vphantom{\displaystyle A^{B^C}}$ &  $x_3 e_3$  \\ & $x_1 e_1 +x_2 e_2$ \; {\rm if} \; $\alpha x_1 ^2 + (\beta+\gamma)x_1x_2 +\delta x_2^2=0$  \\ & $x_1(\delta e_1-\gamma e_2) +x_3 e_3$ {\rm if } \; $\alpha=\beta=0$  \\ & $x_2(-\beta e_1+\alpha e_2) +x_3 e_3$ {\rm if } \; $\gamma=\delta=0$  \\  & 
$x_1 e_1 + x_2 e_2 \pm \sqrt{x_1 ^2 +x_2 ^2} \, e_3$ \; {\rm if} \; $\gamma x_1 ^2 + (\delta-\alpha)x_1x_2 -\beta x_2 ^2=0$\\

\hline  $\mathfrak{g} _6 $ with $(\gamma,\delta)\neq (\beta,\alpha)$ $\vphantom{\displaystyle A^{B^C}}$ &  $x_1 e_1$  \\ & $x_2 e_2 +x_3 e_3$ \; {\rm if} \; $\alpha x_2 ^2 + (\gamma-\beta)x_2x_3 -\delta x_3^2=0$  \\ & $x_1 e_1 +x_2(\delta e_2+\gamma e_3)$ {\rm if } \; $\alpha=\beta=0$  \\ & $x_1 e_1 + x_3(\beta e_2+\alpha e_3)$ {\rm if } \; $\gamma=\delta=0$  \\  & 
$x_1 e_1 + x_2 e_2 \pm \sqrt{x_1 ^2 +x_2 ^2} \, e_3$ \, {\rm if} \; $\gamma x_2 ^2 + (\alpha -\delta)x_2x_3 -\beta x_3 ^2=0$\\

\hline  $\mathfrak{g} _7 $ with $\alpha =0 \neq \gamma$ $\vphantom{\displaystyle A^{B^C}}$ &  $x_1 e_1$  \\  & $x_2 (e_2 + e_3)$ \\ & $-\frac{\delta}{\gamma} (x_2-x_3) e_1 +x_2e_2+ x_3 e_3$ {\rm if } \; $\beta=0$  \\ & $\left(\frac{\beta + \gamma}{\delta} x_2 -  \frac{\beta -\gamma}{\delta} x_3 \right) e_1 +x_2 e_2 + x_3 e_3 \; {\rm if \; (3.23) \; holds }$ \\

\hline  $\mathfrak{g} _7 $ with $\alpha \neq 0 = \gamma$ $\vphantom{\displaystyle A^{B^C}}$ &  $x_2 (e_2 + e_3)$ \\  &  $x_2( e_2 - e_3) \; {\rm if } \; \beta=0$ \\ & 
$x_3 \left\{ 2\beta (\alpha -\delta ) e_1 + [\beta ^2 -(\alpha -\delta )^2] e_2 +
[\beta ^2 +(\alpha -\delta )^2] e_3 \right\}$ \\

\hline
\end{tabular}
\end{center}

\medskip
All information about existence of null and linearly independent homogeneous geodesic through a point, can be easily derived in the different cases from Table III above. It is worthwhile to note that many conditions which appear in Table III, are strictly related to the value of the isomorphism invariant $D=\frac{4(\alpha \delta -\beta \gamma)}{(\alpha +\delta)^2}$. For example, taking into account that, by (3.5), $\alpha \gamma +\beta \delta =0$, we have that, for a Lie algebra $\mathfrak{g} _5 $ with $(\gamma,\delta) \neq (-\beta,\alpha)$, geodesic vectors either of the form $x_1(\delta e_1-\gamma e_2) +x_3 e_3$ or $x_2(-\beta e_1+\alpha e_2) +x_3 e_3$ exist if and only if $D=0$, and geodesic vectors of the form $x_1 e_1 + x_2 e_2 \pm \sqrt{x_1 ^2 +x_2 ^2} \, e_3$ exist if and only if $D \leq 1$.

The most interesting results concerning the existence of linearly independent and null  homogeneous geodesics are summarized in the following Theorems 3.6 and 3.7.

\begin{theorem}
Let $G$ be a three-dimensional non-unimodular Lie group, equipped with a left-invariant   Lorentzian metric $g$, and $\mathfrak{g}$ the Lie algebra of $G$. Assume $(G,g)$ is not symmetric.  

\begin{itemize}
\item When $\mathfrak{g}=\mathfrak{g} _5 $ with $(\gamma,\delta) \neq (-\beta,\alpha)$, 
then through any point of $G$ there are three linearly independent homogeneous geodesics, unless one of the following cases occurs:

\medskip 
(i) $(\alpha,\beta) \neq (0,0) \neq (\gamma,\delta)$, $(\alpha-\delta)^2+4\beta \gamma <0$ and $(\beta+\gamma)^2-4\alpha\delta=0$. In this case, there are two linearly independent homogeneous geodesics.

\medskip 
(ii) $(\alpha,\beta) \neq (0,0) \neq (\gamma,\delta)$, $(\alpha-\delta)^2+4\beta \gamma <0$ and $(\beta+\gamma)^2-4\alpha\delta <0$. In this case, there is just one homogeneous geodesic through a point.
\item When $\mathfrak{g}=\mathfrak{g} _6 $ with $(\gamma,\delta) \neq (\beta,\alpha)$, 
then through any point of $G$ there are three linearly independent homogeneous geodesics, unless:

\medskip 
(i) $(\alpha,\beta) \neq (0,0) \neq (\gamma,\delta)$, $(\alpha-\delta)^2+4\beta \gamma <0$ and $(\beta-\gamma)^2+4\alpha\delta=0$. In this case, there are two linearly independent homogeneous geodesics.

\medskip 
(ii) $(\alpha,\beta) \neq (0,0) \neq (\gamma,\delta)$, $(\alpha-\delta)^2+4\beta \gamma <0$ and $(\beta-\gamma)^2+4\alpha\delta <0$. In this case, there is just one homogeneous geodesic through a point.
\item When $\mathfrak{g}=\mathfrak{g} _7 $ with $\alpha =0 \neq \gamma$, 
then through any point of $G$ there are three linearly independent homogeneous geodesics if and only if either $\beta =0$ or $D \leq 1$. If $\beta \neq 0$ and $D >1$, there are only two linearly independent homogeneous geodesics.
\item When $\mathfrak{g}=\mathfrak{g} _7 $ with $\alpha \neq 0 = \gamma$, 
then through any point of $G$ there are three linearly independent homogeneous geodesics if and only if either $\beta =0$. If $\beta \neq 0$, there are only two linearly independent homogeneous geodesics.

\end{itemize}

\end{theorem}

\begin{theorem}
Consider a three-dimensional non-unimodular Lie group $G$, with Lie algebra $\mathfrak{g}$, equipped with a left-invariant Lorentzian metric $g$. Assume $(G,g)$ is not symmetric.  

\begin{itemize}
\item When $\mathfrak{g}=\mathfrak{g} _5 $, $G$ admits null homogeneous geodesics through a point, except when conditions $(\gamma,\delta) \neq (-\beta,\alpha)$, $(\alpha,\beta) \neq (0,0) \neq (\gamma,\delta)$ and $(\alpha-\delta)^2+4\beta \gamma <0$ are simultaneously satisfied. 
\item When $\mathfrak{g}=\mathfrak{g} _6 $, $G$ admits null homogeneous geodesics through a point, except when conditions $(\gamma,\delta) \neq (\beta,\alpha)$, $(\alpha,\beta) \neq (0,0) \neq (\gamma,\delta)$, $(\alpha-\delta)^2+4\beta \gamma <0$ and $(\beta-\gamma)^2+4\alpha\delta <0$ are simultaneously satisfied. 
\item When $\mathfrak{g}=\mathfrak{g} _7 $, $G$ always admits at least a null homogeneous geodesic through a point. In particular, if $\alpha \neq 0=\gamma$, then all homogeneous geodesics are null.
\end{itemize}

\end{theorem}

\bigskip \noindent
\section{\red Three-dimensional naturally reductive and g.o. Lorentzian spaces}

In the previous Section, we determined all geodesic vectors for the different possible
forms of the Lie algebra of a three-dimensional Lorentzian non-unimodular Lie group. A corresponding investigation was made in [CM2] in the unimodular case. Because of Theorem 1.1, every (connected, simply connected) three-dimensional homogeneous Lorentzian manifold $M$ is isometric to a Lie group equipped with a left-invariant Lorentzian metric, unless $M$ is symmetric. If $M$ is symmetric, then in particular it is naturally reductive and a g.o. space. Hence, we actually know the sets of homogeneous geodesics through a point for {\em all} three-dimensional homogeneous Lorentzian manifolds, and this allows us to classify three-dimensional Lorentzian g.o. and naturally reductive spaces. 

The results of Section 4 of [CM2], as concerns the existence of g.o. spaces, can be summarized in the following

\begin{theorem}{\bf [CM2]} 
A three-dimensional unimodular Lorentzian Lie group $(G,g)$ is a 
g.o. space if and only if either it is symmetric, or its Lie algebra $\mathfrak{g}$ is one of the following:
\begin{itemize}
\item $\mathfrak{g}=\mathfrak{g} _3$, with either $\alpha = \beta \neq \gamma$, $\alpha = \gamma \neq \beta$ or $\beta = \gamma \neq \alpha$.
\item $\mathfrak{g}=\mathfrak{g} _4 $, with $\alpha = \beta - \varepsilon$.
\end{itemize}
Moreover, in all these cases, $G$ is also naturally reductive. 
\end{theorem}

On the other hand, Section 3 here permits to conclude that a three-dimensional non-unimodular Lorentzian Lie group is never a g.o. space, unless it is symmetric. (The same result is true in the Riemannian case [TV].) Therefore, the following classification result holds:

\begin{theorem} 
Let $(M,g)$ be a connected, simply connected three-dimensional Lorentzian manifold. The following properties are equivalent:

\medskip (i)  $(M,g$) is a g.o. space.

\medskip (ii)  $(M,g$) is naturally reductive.

\medskip (iii) Either $(M,g$) is symmetric, or it is isometric to a unimodular Lie group $G$, equipped with a left-invariant Lorentzian metric, having one of the following Lie algebras:
\begin{itemize}
\item $\mathfrak{g}=\mathfrak{g} _3$, with either $\alpha = \beta \neq \gamma$, $\alpha = \gamma \neq \beta$ or $\beta = \gamma \neq \alpha$.
\item $\mathfrak{g}=\mathfrak{g} _4 $, with $\alpha = \beta - \varepsilon$.
\end{itemize}
\end{theorem}

Unimodular Lie groups admitting one of the Lie algebras listed in Theorem 4.2, can be easily deduced from Tables I and II. Hence, we can get the explicit classification of three-dimensional non-symmetric naturally reductive Lorentzian spaces.

\begin{theorem} 
A three-dimensional connected, simply connected Lorentzian manifold $(M,g)$ is a non-symmetric naturally reductive space if and only if it is isometric to 

\medskip
a) $SL(2,\Real)$,

\medskip
b) $SU(2)$,

\medskip
c) $H_3$,

\medskip\noindent
equipped with a suitable left-invariant Lorentzian metric.
\end{theorem}

It is worthwhile to compare this result with Theorem 6.5 in [TV], where it was proved that $SL(2,\Real)$, $SU(2)$ and $H_3$, equipped with a suitable left-invariant Riemannian metric, are the only three-dimensional naturally reductive non-symmetric  Riemannian spaces.

\bigskip

\begin{center}
\textsc{Dipartimento di Matematica "E. De Giorgi",\\
Universit\`{a} degli Studi di Lecce,\\ Via Provinciale Lecce-Arnesano, \\
73100 Lecce, ITALY.} \\
\textit{E-mail addresses}: giovanni.calvaruso@unile.it, rosanna.marinosci@unile.it
\end{center}

\end{document}